\documentclass[final]{amsart}

\usepackage{stmaryrd}
\usepackage{amssymb}
\usepackage{amsmath}
\usepackage{amsthm}
\usepackage{amsfonts}
\usepackage{amscd}
\usepackage{pifont}
\usepackage{boxedminipage}
\usepackage{ifthen}
\usepackage[usenames]{color}
\usepackage{graphicx}

\newtheorem{Lemma}{Lemma}[section]
\newtheorem{Theorem}[Lemma]{Theorem}

\newtheorem{Proposition}[Lemma]{Proposition}

\theoremstyle{definition}

\newtheorem{Remark}[Lemma]{Remark}

\numberwithin{equation}{section}

\newcommand{\supp}{\operatorname{supp}}

\newcommand{\Fix}[1][J]{\ensuremath{{\mathfrak c}_{#1}}}

\newcommand{\rad}{\operatorname{rad}}

\newcommand{\inv}{^{-1}}

\newcommand{\mc}{\mathcal}
\newcommand{\mb}{\mathbb}
\newcommand{\bs}{\boldsymbol}
\newcommand{\Arrangement}{\mc A}

\newcommand{\DescentAlgebra}[1][W]{\ensuremath{\Sigma_k(#1)}}
\newcommand{\Faces}{\mc F}
\newcommand{\InvariantSubalgebra}[1][W]{\ensuremath{(k\Faces)^{#1}}}
\newcommand{\InvariantQuiverAlgebra}[1][W]{\ensuremath{(k\Quiver)^{#1}}}
\newcommand{\InvariantQuiver}[1][W]{\ensuremath{\Quiver_{#1}}}
\newcommand{\IntersectionLattice}{\mc L}
\newcommand{\Orbit}{\mc O}
\newcommand{\Quiver}{\mc Q}
\newcommand{\CompleteSystemPunctuated}[1][c]{{#1}omplete system of primitive orthogonal idempotents}
\newcommand{\CompleteSystem}[1][c]{{\CompleteSystemPunctuated[#1] }}

\renewcommand{\eqref}[1]{(\ref{#1})}
\newcommand{\LoewyLength}[1][W]{\operatorname{LL}(#1)}

\newcommand{\DescentAlgebraD}[1][2m+1]{\ensuremath{\DescentAlgebra[D_{#1}]}}
\newcommand{\InvariantQuiverD}[1][2m+1]{\ensuremath{\InvariantQuiver[D_{#1}]}}
\newcommand{\InvariantSubalgebraD}[1][2m+1]{\ensuremath{\InvariantSubalgebra[D_{#1}]}}

\newcommand{\adots}{\mathord{\cdot}\mathord{\cdot}\mathord{\cdot}}
\newcommand{\arrow}{\mathord{\rightarrow}}
\newcommand{\arrows}{\arrow\adots\arrow}
\newcommand{\Path}[3][1]{
\ifthenelse{\equal{#1}{2}}%
 {\ensuremath{\left({#2}_0\arrow{#2}_1\arrows{#2}_{#3}\right)}}%
 {\ensuremath{\left({#2}_0\arrows{#2}_{#3}\right)}}%
}

\newcommand{\SetPartitionD}%
  {\ensuremath{\{B_1,\ldots,B_r;C\}}}

\newcommand{\Even}[1]{\ensuremath{\operatorname{Even}({#1})}}
\newcommand{\Odd}[1]{\ensuremath{\operatorname{Odd}({#1})}}

\newcommand{\dimDiff}[1][j]{\ensuremath{\dim(X_{#1-1})-\dim(X_{#1})}}

\newenvironment{note}[1][Note]
 {\begin{center}\begin{boxedminipage}{4.5in}\setlength{\parindent}{1em}\noindent\textbf{#1. }}
 {\end{boxedminipage}\end{center}}

\newenvironment{Notation}{\noindent\textsc{Notation.}}{\smallskip}

\newcommand{\defn}[1]{\textit{\textbf{#1}}}

\allowdisplaybreaks

\begin{document}

\title[The Loewy length of $\DescentAlgebraD$]{The Loewy length of
the\\ descent algebra of $D_{2m+1}$.}
\author{Franco V Saliola}
\address{Laboratoire de Combinatoire et d'Informatique Math\'ematique\\
	Universit\'e du Qu\'ebec \`a Montr\'eal\\
	Case Postale 8888, succursale Centre-ville\\
	Montr\'eal, Qu\'ebec H3C 3P8\\
	Canada}
\email{saliola@gmail.com}

\begin{abstract}
In this article the Loewy length of the descent algebra of $D_{2m+1}$
is shown to be $m+2$, for $m \geq 2$, by providing an upper bound that
agrees with the lower bound in \cite{BonnafePfeiffer2006}. 

The bound is obtained by showing that the length of the longest path
in the quiver of the descent algebra of $D_{2m+1}$ is at most $m+1$.
To achieve this bound, the geometric approach to the descent algebra
is used, in which the descent algebra of a finite Coxeter group $W$ is
identified with an algebra associated to the reflection arrangement of
$W$.
\end{abstract}

\maketitle


\section{Introduction}

C\'edric Bonnaf\'e and G\"otz Pfeiffer determined the Loewy length of
the descent algebra $\DescentAlgebra$ for irreducible finite Coxeter
groups $W$ of all types except $D_{2m+1}$ \cite{BonnafePfeiffer2006}.
(The case of the symmetric group was established earlier; see
\cite[Theorems 5.6 and 5.7]{GarsiaReutenauer1989} and
\cite[Theorem 3.4]{Atkinson1992}).
For type $D_{2m+1}$, Bonnaf\'e and Pfeiffer prove that the Loewy
length of $\DescentAlgebraD$ is at least $m+2$, and state that they
suspect this is an equality. In this paper, we show that $m+2$ is also
an upper bound---by showing that the length of the longest path in
the quiver of $\DescentAlgebraD$ is at most $m+1$---confirming their
suspicion. 

We briefly outline the argument, and the structure of the paper. The
plan is to use the \emph{geometric approach} to the descent algebra:
the descent algebra $\DescentAlgebra$ can be identified with a
subalgebra of an algebra $k\Faces$ associated to the reflection
arrangement of $W$. This is explained in Section
\ref{s:RefArrs+DesAlg}. Specifically, there is an action of $W$ on
$k\Faces$ and an anti-isomorphism between $\DescentAlgebra$ and the
$W$-invariant subalgebra $\InvariantSubalgebra$. After defining
quivers and path algebras in Section \ref{s:QuiversAndPathAlgs},
Section \ref{s:EquivariantMap} describes the quiver $\Quiver$ of
$k\Faces$ and the construction of a $W$-equivariant surjection
$\varphi: k\Quiver\to k\Faces$. This results in a surjection
$(k\Quiver)^W\twoheadrightarrow\InvariantSubalgebra$ that we use to
gain information about the quiver $\InvariantQuiver$ of
$\InvariantSubalgebra$ in Section \ref{s:QuiverOfInvSub}. We then
specialize in Section \ref{s:LoewyLengthTypeD} to the irreducible
finite Coxeter group of type $D$ and bound the length of the longest
path in the quiver $\InvariantQuiverD$. The reader familiar with the
above theory may decide to begin in Section \ref{s:LoewyLengthTypeD}.

Throughout this paper $k$ denotes a field whose characteristic is zero
or does not divide the order of the Coxeter group. 

\newpage

\section{The Geometric Approach to the Descent Algebra}
\label{s:RefArrs+DesAlg}

We begin by recalling the definition of Coxeter systems and the
descent algebra. We then explain the connection between the descent
algebra and the face semigroup algebra of the reflection arrangement
of the Coxeter group.

\subsection{Coxeter systems and reflection arrangements}
\label{ss:CoxSys+RefArrs}

Let $V$ be a finite dimensional real vector space. A \defn{finite
Coxeter group} $W$ is a finite group generated by a set of reflections
of $V$. The \defn{reflection arrangement} of $W$ is the hyperplane
arrangement $\Arrangement$ consisting of the hyperplanes of $V$ fixed
by some reflection in $W$. 

Let $W$ denote a finite Coxeter group with reflection arrangement
$\Arrangement$ and let $c$ denote a connected component of the
complement of $\bigcup_{H\in\Arrangement} H$ in $V$. A \defn{wall} of
$c$ is a hyperplane $H\in\Arrangement$ such that $H \cap \overline c$
spans $H$, where $\overline c$ is the closure of $c$ in $V$. Let $S
\subseteq W$ denote the set of reflections in the walls of $c$. Then
$S$ is a generating set of $W$ \cite[{\S}I.5A]{Brown1989} and the pair
$(W,S)$ is called a \defn{Coxeter system} with \defn{fundamental
chamber} $c$.

\subsection{The descent algebra}
\label{ss:DesAlg}

Fix a Coxeter system $(W,S)$. For $J \subseteq S$, let $W_J = \langle
J \rangle$ denote the subgroup of $W$ generated by the elements in
$J$. Each coset of $W_J$ in $W$ contains a unique element of minimal
length, where the \defn{length} $\ell(w)$ of an element $w$ of $W$ is
the smallest number of generators $s_1, \ldots, s_i \in S$ such that
$w = s_1 \cdots s_i$ \cite[Proposition 1.10(c)]{Humphreys1990}. Let
$X_J$ denote the set of \defn{minimal length coset representatives} of
$W_J$ and let $x_J = \sum_{w \in X_J} w$ denote the formal sum of the
elements of $X_J$. Then $x_J$ is an
element of the group algebra $kW$ of $W$ with coefficients in a field
$k$. Louis Solomon proved that the elements $x_J$, one for each $J
\subseteq S$, form a basis of a subalgebra of $kW$ \cite[Theorem
1]{Solomon1976}. This subalgebra is denoted by $\DescentAlgebra$ and
is called the \defn{descent algebra} of $W$.

\subsection{The faces of $\Arrangement$}
\label{ss:Faces}
For each hyperplane $H \in \Arrangement$, let $H^+$ and $H^-$ denote
the two open half spaces of $V$ determined by $H$. The choice of
labels $H^+$ and $H^-$ is arbitrary, but
fixed throughout. For convenience, let $H^0 = H$. A \defn{face}
of the arrangement $\Arrangement$ is a non-empty intersection of the
form $x = \bigcap_{H \in \Arrangement} H^{\sigma_H(x)}$, where
$\sigma_H(x) \in \{+,0,-\}$ for each hyperplane $H \in \Arrangement$.
The sequence $\sigma(x) = (\sigma_H(x))_{H \in \Arrangement}$ is
called the \defn{sign sequence} of $x$. The set $\Faces$ of all faces
of $\Arrangement$ is a partially order set with partial order given by
$x \leq y$ iff $x \subseteq \overline{y}$, where $\overline{y}$
denotes the closure of the set $y$. A \defn{chamber} of $\Arrangement$
is a face that is maximal with respect to this order. 

\subsection{The intersection lattice}
\label{ss:IntLattice}
For each face $x \in \Faces$, the \defn{support} $\supp(x)$ of $x$ is
the intersection of all hyperplanes in $\Arrangement$ that contain
$x$. Equivalently, $\supp(x)$ is the subspace of $V$ spanned by $x$.
The \defn{dimension} of $x$ is the dimension of the subspace
$\supp(x)$. The \defn{intersection lattice} $\IntersectionLattice$
of $\Arrangement$ is the image of $\supp$; that is,
$\IntersectionLattice = \supp(\Faces)$. The elements
of $\IntersectionLattice$ are subspaces of $V$ and are ordered by
inclusion. With this partial order, $\IntersectionLattice$ is a finite
lattice, where the meet of two subspaces is their intersection, and
the join of two subspaces is the smallest subspace that contains both.
It follows that $\supp: \Faces \to \IntersectionLattice$ is an
order-preserving surjection of posets. (N.B. Some authors order
$\IntersectionLattice$ by reverse inclusion rather than inclusion.)

\subsection{The face semigroup algebra}
\label{ss:FaceSemigroupAlg}
Define the \defn{product of two faces} $x, y \in
\Faces$ to be the face $xy$ with sign sequence $(\sigma_H(xy))_{H \in
\Arrangement}$ given by
\begin{gather*}
\sigma_H(xy) = 
\begin{cases}
\sigma_H(x), & \text{if } \sigma_H(x) \neq 0, \\
\sigma_H(y), & \text{if } \sigma_H(x) = 0,
\end{cases}
\end{gather*}
where $\sigma(x)$ and $\sigma(y)$ are the sign sequences of $x$ and
$y$. This product has a geometric interpretation: the product $xy$ of
two faces $x$ and $y$ is the face entered by moving a small positive
distance along a straight line from a point in $x$ to a point in $y$.
It is straightfoward to verify that this product gives $\Faces$ the
structure of an associative semigroup with identity, and that $x^2 =
x$ and $xyx = xy$ for all $x,y \in \Faces$. (A semigroup satisfying
these identities is called a \emph{left regular band}.)

The semigroup algebra $k\Faces$ is called the \defn{face semigroup
algebra} of $\Arrangement$ over the field $k$. It consists of finite
$k$-linear combinations of elements of $\Faces$ with multiplication
induced by the product defined on elements of $\Faces$. 

\subsection{The invariant subalgebra}
\label{ss:InvSubAlg}

Since $W$ is a group of orthogonal transformations of the vector space
$V$, there is an action of $W$ on $V$ defined by setting $w(\vec v)$
to be the image of $\vec v \in V$ under the transformation $w$. Under
this action the set $\Arrangement$ is permuted \cite[Proposition
1.2]{Humphreys1990}, so there is an induced action of $W$ on
$\IntersectionLattice$ and on $\Faces$. The action preserves the
semigroup structure of $\Faces$, so it extends linearly to an action
on $k\Faces$. Let $\InvariantSubalgebra$ denote the subalgebra of
$k\Faces$ consisting of the elements of $k\Faces$ fixed by all
elements of $W$: 
\begin{gather*}
\InvariantSubalgebra = \Big\{ a \in k\Faces : w(a) = a \text{ for all
} w\in W \Big\}.
\end{gather*}

The following was first proved by T. P. Bidigare \cite{Bidigare1997}.
Another proof was given by K. S. Brown and can be found in
\cite[Theorem 7]{Brown2000} or \cite[Theorem 2.7]{Saliola2007a}.

\begin{Theorem}
\label{t:Bidigare}
Let $W$ be a finite reflection group and let $k\Faces$ denote the face
semigroup algebra of the reflection arrangement of $W$. The
$W$-invariant subalgebra $\InvariantSubalgebra$ is anti-isomorphic to
the descent algebra $\DescentAlgebra$ of $W$.
\end{Theorem}

We briefly describe an anti-isomorphism. The faces of the fundamental
chamber $c$ are parametrized by the subsets of $S$: if $J \subseteq
S$, then there is a unique face $\Fix[J]$ of the fundamental chamber
$c$ that is fixed by all elements of $J$ \cite[{\S}I.5F]{Brown1989}.
Moreover, every face of $\Arrangement$ is in the $W$-orbit of a unique
face of $c$ \cite[{\S}I.5F]{Brown1989}. So if $\Orbit_J$ denotes the
$W$-orbit of $\Fix[J]$, then the elements $\bs x_J = \sum_{y \in
\Orbit_J} y$, one for each $J \subseteq S$, form a basis of
$\InvariantSubalgebra$. The function $\InvariantSubalgebra \to
\DescentAlgebra$ defined by mapping $\bs x_J$ to $x_J$ is an
anti-isomorphism.

\section{Quivers and Path Algebras}
\label{s:QuiversAndPathAlgs}

Let $k$ be a field and $A$ a finite dimensional $k$-algebra. 

\subsection{\CompleteSystemPunctuated[C]}
An element $a \in A$ is an \defn{idempotent} if $e^2 = e$. Two
idempotents $e,f\in A$ are \defn{orthogonal} if $ef = 0 = fe$. An
idempotent $e\in A$ is \defn{primitive} if it cannot be written as $e
= f + g$ with $f$ and $g$ non-zero orthogonal idempotents of $A$. A
\defn{\CompleteSystemPunctuated} of $A$  is a set $\{e_1, e_2, \ldots,
e_n\}$ of primitive idempotents of $A$ that are pairwise orthogonal
and that sum to $1_A$.

\subsection{The quiver of a split basic algebra}
\label{ss:QuiverDefn}

The \defn{Jacobson radical} $\rad A$ of $A$ is the smallest ideal of
$A$ such that $A/\rad A$ is semisimple. If $A/\rad A$ is a direct
product of copies of $k$, then $A$ is a \defn{split basic
$k$-algebra}. Equivalently, $A$ is a split basic algebra if and only
if all the simple $A$-modules are of dimension one. 

The \defn{quiver} $Q$ of a split basic finite dimensional $k$-algebra
$A$ is the finite directed graph constructed as follows. Let $\{e_v: v
\in \mc I\}$ be a \CompleteSystem of $A$, where $\mc I$ is some index
set. The vertex set of $Q$ is the index set $\mc I$, so there is one
vertex $v$ in $Q$ for each idempotent $e_v$. If $x,y \in \mc I$, then
the number of arrows $x\arrow y$ is $\dim_k
e_y\big(\rad(A)/\rad^2(A)\big)e_x$. This construction does not depend
on the \CompleteSystemPunctuated, so $Q$ is canonically defined.

\subsection{The path algebra}
\label{ss:PathAlgebraDefn}
The \defn{path algebra} $kQ$ of a quiver $Q$ is the $k$-algebra with
basis the set of paths in $Q$ and with multiplication defined on
paths by 
\begin{gather*}
\Path ws \cdot \Path vr =
\begin{cases}
(v_0 \arrow\cdots\arrow v_r \arrow w_1 \arrow\cdots\arrow w_s),
& \text{if }w_0 = v_r, \\
\hfill0\hfill,
& \text{if }w_0 \neq v_r,
\end{cases}
\end{gather*}
where $\Path ws$ and $\Path vr$ are paths in $Q$. If $F$ denotes the
ideal of $kQ$ generated by the arrows of $Q$, then an ideal $I
\subseteq kQ$ is said to be \defn{admissible} if $F^m \subseteq I
\subseteq F^2$ for some $m\in\mb N$.

If $Q$ is the quiver of $A$, then there is a surjection $\varphi: kQ
\to A$ defined by mapping each vertex $x$ to the idempotent $e_x$ and
by mapping the arrows from $x$ to $y$ to elements in
$e_y \rad(A) e_x$ whose image in $e_y\big(\rad(A)/\rad^2(A)\big)e_x$
forms a basis of the quotient space. Moreover, $\ker \varphi$ is an
admissible ideal of $kQ$.


\subsection{Loewy length} 
The \defn{Loewy length} $\LoewyLength[A]$ of a finite dimensional
$k$-algebra $A$ is the smallest $l\in\mb N$ such $(\rad A)^l = 0$.
The following observation is pertinent.

\begin{Lemma}
\label{l:LoewyLength}
Suppose $Q$ is a finite acyclic quiver. If $A \cong kQ/I$ for some
quiver $Q$ and some admissible ideal $I$ of $kQ$, then
$\LoewyLength[A]\leq l+1$, where $l$ is the length of the longest path
in $Q$. 
\end{Lemma}

\begin{proof}
Let $\varphi: kQ \to A$ denote a surjection with kernel $I$. If $F
\subseteq kQ$ denotes the ideal generated by the arrows of $Q$, then
$\varphi(F^l) = (\rad A)^l$ for all $l\geq1$ \cite[Corollary
2.11]{Assem2006}. So if $l\in\mb N$ is the length of the longest path
in $Q$, then $F^{l+1} = 0$. Hence, $(\rad A)^{l+1} = 0$. Thus,
$\LoewyLength[A]\leq l+1$.
\end{proof}

\section{A $W$-Equivariant Surjection}
\label{s:EquivariantMap}

\begin{Notation}
Throughout this section: $(W,S)$ is a Coxeter system with fundamental
domain $c$; $\Arrangement$ is the reflection arrangement of $W$;
$\IntersectionLattice$ is the intersection lattice of $\Arrangement$;
and $k\Faces$ is the face semigroup algebra of $\Arrangement$,
where $k$ is a field whose characteristic does not divide the order
of $W$.
\end{Notation}

In this section we recall the construction a \CompleteSystem for
$\InvariantSubalgebra$ and the construction of a $W$-equivariant
surjection $\varphi:k\Quiver\to k\Faces$, where $\Quiver$ is the
quiver of $k\Faces$. 

\subsection{The orbit poset}
\label{ss:OrbitPoset}
For each $x \in \Faces$ let $\Orbit_x = \{w(x) : w \in W\}$ denote the
$W$-orbit of $x$, and for each $X \in \IntersectionLattice$ let
$\Orbit_X = \{ w(X) : w \in W \}$ denote the $W$-orbit of $X$. The
$W$-orbits of elements of $\IntersectionLattice$ form a poset
$\IntersectionLattice/W = \{ \Orbit_X : X \in \IntersectionLattice \}$
with partial order given by $\Orbit_X \leq \Orbit_Y$ if and only if
there exists $w \in W$ with $w(X) \leq Y$.

\begin{Remark}
\label{r:IntLatAndSubsetsOfS}
The poset $\IntersectionLattice/W$ is isomorphic to a poset of
equivalence classes of subsets of $S$. Indeed, define a relation on
subsets $J,K \subseteq S$ by setting $J \sim K$ if and only if
$\supp(\Fix[J])$ and $\supp(\Fix[K])$ belong to the same orbit.
Equivalently, $J \sim K$ if and only if $W_J$ and $W_K$ are conjugate
subgroups of $W$. The poset $S/\mathord{\sim}$, with partial order
induced by reverse inclusion of subsets of $S$, is isomorphic to
$\IntersectionLattice/W$. 
\end{Remark}

\subsection{\CompleteSystemPunctuated[C]}
\label{ss:ComSys}
The construction requires, for each $X\in\IntersectionLattice$, a
linear combination $\ell(X)$ of faces of support $X$ with coefficients
summing to $1$. Moreover, the elements $\ell(X)$ need to satisfy
$w(\ell(X)) =
\ell(w(X))$ for all $w\in W$. 

We provide one example of such elements; see \S3.4 of
\cite{Saliola2007a} for other examples.
For every orbit $\Orbit\in\IntersectionLattice/W$, fix a face
$f_\Orbit$ such that $\supp(f_\Orbit)\in\Orbit$. 
For each $X\in\IntersectionLattice$, let $f_X = f_{\Orbit_X}$
and define
\begin{gather}
\label{e:lambda}
\ell(X) = 
\frac1{\lambda_X}
 \sum_{ z \in \Orbit_{f_X} \atop \supp(z) = X } z, 
\quad \text{ where }
\lambda_X 
=|\{ z \in \Orbit_{f_X} : \supp(z) = X \}|.
\end{gather}
Note that $\lambda_X$ is the index of the stabilizer subgroup $W_x =
\{w\in W: w(x)=x\}$ of $x$ in the stabilizer subgroup $W_X=\{w\in
W:w(X)=X\}$ of $X$, where $x$ is any face with support $X$. Hence,
$\lambda_X$ depends only on the orbit of $X$ and so the elements
$\ell(X)$
satisfy $w(\ell(X)) = \ell(w(X))$ for all $w\in W$.

Define elements $e_X\in k\Faces$, one for each $X \in
\IntersectionLattice$, recursively by the formula 
\begin{gather}
\label{e:ComSysForFaceSemiAlg}
e_X = \ell(X) - \ell(X) \sum_{Y > X} e_Y.
\end{gather}
These elements form a \CompleteSystem for $k\Faces$
\cite[Theorem 5.2]{Saliola2006a}. Moreover, they
satisfy $w(e_X) = e_{w(X)}$ for all $w\in W$ and all
$X\in\IntersectionLattice$, so the elements
\begin{gather}
\label{e:ComSysForInvSub}
\varepsilon_\Orbit = \sum_{X\in\Orbit} e_X,
\end{gather}
one for each $\Orbit\in\IntersectionLattice/W$, form a \CompleteSystem
for $\InvariantSubalgebra$ \cite[Theorem 3.7]{Saliola2007a}.

\begin{Remark}
The above leads to a construction of a \CompleteSystem directly within
the descent algebra \DescentAlgebra. Let $S/\mathord{\sim}$ denote the
poset defined in Remark \ref{r:IntLatAndSubsetsOfS}. For each $\Orbit
\in S/\mathord{\sim}$, fix a subset $J_\Orbit\subseteq S$ with
$J_\Orbit \in \Orbit$ and define elements $\varepsilon_\Orbit$, one
for each $\Orbit\in S/\mathord{\sim}$, recursively by the formula
\begin{gather*}
\varepsilon_\Orbit 
  = \frac1{\lambda_\Orbit} x_{J_\Orbit} 
  - \sum_{\Orbit'>\Orbit} \varepsilon_{\Orbit'}
     \left(\frac1{\lambda_\Orbit} 
     x_{J_\Orbit} \right),
\end{gather*}
where $\lambda_\Orbit$ is the index of $W_J$ in the normalizer of
$W_J$. These elements correspond, under the anti-isomorphism of
\S\ref{ss:InvSubAlg}, to the elements defined in Equation
\ref{e:ComSysForInvSub} for a suitable choice of $f_\Orbit$. 
Therefore, they form a \CompleteSystem for \DescentAlgebra.
See \cite[Proposition 3.9]{Saliola2007a} for details.
\end{Remark}

\subsection{The quiver of $k\Faces$}
\label{ss:QuiverOfFaceSemiAlg}

The quiver of $k\Faces$ is the directed graph $\Quiver$ constructed as
follows. The vertex set of $\Quiver$ is $\IntersectionLattice$, and
there is exactly one arrow $X\arrow Y$, for
$X,Y\in\IntersectionLattice$, if and only if $Y \lessdot X$. There
exists a surjection $\varphi: k\Quiver \to k\Faces$ with kernel
generated by the sum of all the paths in $\Quiver$ of length two. A
proof of this for any central hyperplane arrangement can by found in
\cite{Saliola2006a}. Below we recall the construction
of a $W$-equivariant $\varphi: k\Quiver\twoheadrightarrow k\Faces$.
See \cite[{\S}4]{Saliola2007a} for details.

Define an action of $W$ on the path algebra $k\Quiver$ as follows.
Fix an orientation $\epsilon_X$ on each subspace $X \in
\IntersectionLattice$. Thus, $\epsilon_X$ is a map that assigns $1$ or
$-1$ to a basis of $X$ depending on whether the basis is positively or
negatively oriented. For $w \in W$ and $X \in \IntersectionLattice$, let 
\begin{gather*}
\sigma_X(w) = \epsilon_X\Big(\vec x_1, \ldots, \vec x_s\Big)
 \epsilon_{w(X)}\Big(w(\vec x_1), \ldots, w(\vec x_s)\Big),
\end{gather*}
where $\vec x_1, \ldots, \vec x_s$ is a basis of the subspace $X$.
Note that if $w(X) = X$, then $\sigma_X(w)$ is $1$ if the restriction
$w|_X$ of $w$ to $X$ is orientation-preserving, and is $-1$ otherwise.
For $w \in W$ and a path $\Path Xt$ in $\Quiver$ define
\begin{align*}
w \Path Xt
= \sigma_{X_0}(w)\sigma_{X_t}(w) \Big(w(X_0) \arrows w(X_t)\Big),
\end{align*}
where $w(X_i)$ is the image of $X_i\in\IntersectionLattice$ under
the action of $W$ on $\IntersectionLattice$. 
The following is Theorem 4.7 of \cite{Saliola2007a}.

\begin{Theorem}
\label{t:DefnOfPhi}
There exists a $W$-equivariant $k$-algebra surjection $\varphi:
k\Quiver \twoheadrightarrow k\Faces$. That is, $\varphi$ satisfies
$\varphi(w(a)) = w(\varphi(a))$ for all $w \in W$ and all $a \in
k\Quiver$. In addition, $\ker(\varphi)$ is generated by the sum of all
the paths of length two in $\Quiver$.
\end{Theorem}

We briefly recall that construction. Let $\epsilon_X$
denote the orientations chosen above. Define $\varphi$ on each vertex
$X$ and arrow $X\arrow Y$ of $\Quiver$ by
\begin{align*}
\varphi(X) = e_X\qquad\text{ and }\qquad
\varphi(X\arrow Y) 
 = \lambda_Y \, e_Y\!\left(\sum_{x \gtrdot y} [y:x] x \right)\!e_X,
\end{align*}
where $\lambda_Y$ is defined in Equation \eqref{e:lambda}, 
where $y$ is any face of support $Y$ and
\begin{align*}
[y:x] = \epsilon_{\supp(y)}(\vec y_1, \ldots, \vec y_t)
 \epsilon_{\supp(x)}(\vec y_1, \ldots, \vec y_t, \vec x_1),
\end{align*}
where $\vec y_1, \ldots, \vec y_t$ is a basis of $\supp(y)$ and $\vec
x_1$ is a vector in $x$. Then $\varphi$ extends linearly and
multiplicatively to a $W$-equivariant $k$-algebra surjection $\varphi:
k\Quiver \to k\Faces$. 

\section{On the Quiver of \InvariantSubalgebra}
\label{s:QuiverOfInvSub}

We continue with the notation of the previous section and let
$\InvariantQuiver$ denote the quiver of $\InvariantSubalgebra$. The
quiver $\InvariantQuiver$ is not known for arbitrary $W$, but the
idempotents of Equation \eqref{e:ComSysForInvSub} and the surjection
of Theorem \ref{t:DefnOfPhi} provide some information about the
structure of $\InvariantQuiver$. In the next section we specialize to
$W$ of type $D$.

The following result is our main tool.

\begin{Lemma}
\label{l:ConditionForNoArrow}
If for every path $P$ in $\Quiver$ that begins at a vertex in
$\Orbit'\in\IntersectionLattice/W$ and ends at a vertex in
$\Orbit\in\IntersectionLattice/W$ there exists $w \in W$ such that
$w(P)=-P$, then there is no arrow from $\Orbit'$ to $\Orbit$ in
$\InvariantQuiver$.
\end{Lemma}

\begin{proof}
If there is an arrow $\Orbit'\arrow\Orbit$, then the vector space
$\varepsilon_\Orbit\InvariantSubalgebra\varepsilon_{\Orbit'}$ is
nonzero (see \S\ref{ss:QuiverDefn}). We'll show that this vector space
is zero if the hypothesis holds.

For each $\Orbit\in\IntersectionLattice/W$, let $\nu_\Orbit =
\sum_{X\in\Orbit} X \in k\Quiver$. 
Let $\varphi: k\Quiver \to k\Faces$ denote the $W$-equivariant
surjection of Theorem \ref{t:DefnOfPhi}.
Then $\varphi$ restricts to a surjection 
\begin{align*}
\nu_\Orbit \InvariantQuiverAlgebra \nu_{\Orbit'}
\twoheadrightarrow
\varepsilon_\Orbit \InvariantSubalgebra \varepsilon_{\Orbit'}.
\end{align*}
We'll show that $\nu_\Orbit\InvariantQuiverAlgebra\nu_{\Orbit'}=0$.
This subspace is spanned by elements of the form $\sum_{P'\in\Orbit_P}
P'$, where $P$ is a path of $\Quiver$ that begins at a vertex in
$\Orbit'$ and ends at a vertex in $\Orbit$, and where $\Orbit_P$ is
the $W$-orbit of $P$. The hypothesis implies $w(P)=-P$ for some $w\in
W$, so
\begin{gather*}
\sum_{P'\in\Orbit_P} P' 
= \sum_{P'\in\Orbit_P} w(P')
= \sum_{P'\in\Orbit_{-P}} P'
= -\sum_{P'\in\Orbit_P} P'.
\end{gather*}
Therefore, $\sum_{P'\in\Orbit_P} P'=0$. So
$\nu_\Orbit\InvariantQuiverAlgebra\nu_{\Orbit'}=0$.
\end{proof}

Our first result on the structure of $\InvariantQuiver$ 
shows that it contains no oriented cycles.
\begin{Proposition}
\label{p:NoOrientedCycles}
There is exactly one vertex in $\InvariantQuiver$ for each element of
$\IntersectionLattice/W$. If $\Orbit' \arrow \Orbit$ is an arrow in
$\InvariantQuiver$, then $\Orbit \leq \Orbit'$ in
$\IntersectionLattice/W$. In particular, $\InvariantQuiver$ does not
contain any oriented cycles.
\end{Proposition}
\begin{proof}
Since the elements in Equation \eqref{e:ComSysForInvSub} form a
\CompleteSystem for \InvariantSubalgebra, the vertex set of
$\InvariantQuiver$ is the poset $\IntersectionLattice/W$.

If $\Path Xl$ is a path in $\Quiver$, then $X_l\leq X_0$.
In particular, $\Orbit_{X_l} \leq \Orbit_{X_0}$. So if
$\Orbit\nleq\Orbit'$, then the condition of Lemma
\ref{l:ConditionForNoArrow} is satisfied. Therefore, there is no arrow
from $\Orbit'$ to $\Orbit$ in $\InvariantQuiver$. It follows that
$\InvariantQuiver$ cannot contain an oriented cycle.
\end{proof}


Our next result shows that the quiver $\InvariantQuiver$ contains at
least one isolated vertex.

\begin{Proposition}
\label{p:NoArrowsFromTop}
There are no arrows in $\InvariantQuiver$ beginning at $\{V\}$.
\end{Proposition}
\begin{proof}
Let $\Path Xl$ be a path in $\Quiver$ with $X_0=V$. Let $w\in W$ denote
the reflection in the hyperplane $X_1$. Then 
\begin{align*}
w\Path Xl 
&= \sigma_{X_0}(w)\sigma_{X_l}(w)
\left( w(X_0)\arrow\adots\arrow w(X_l) \right) \\
&= -\Path{X}l.
\end{align*}
By Lemma \ref{l:ConditionForNoArrow}, there is no arrow
in $\InvariantQuiver$ beginning at $\{V\}$.
%
\end{proof}

\section{The Loewy Length of $\DescentAlgebraD$}
\label{s:LoewyLengthTypeD}

\begin{Notation}
Throughout this section: $D_n$ is a Coxeter group of type $D$,
$\Arrangement$ is the reflection arrangement of $D_n$;
$\IntersectionLattice$ is the intersection lattice of $\Arrangement$;
and $k\Faces$ is the face semigroup algebra of $\Arrangement$, where
$k$ is a field whose characteristic does not divide the order of
$D_n$.
\end{Notation}

\subsection{Coxeter groups of type \emph{D}}
\label{ss:TypeDCoxGrp}

For $n\in\mb N$, let $[n] = \{1,2,\ldots,n\}$ and let
$[\pm n]=\{1,2,\ldots,n\} \cup
\{-1,-2,\ldots,-n\}$. A \defn{signed permutation} of $[\pm n]$ is a
permutation $w$ of $[\pm n]$ satisfying $w(-i)=-w(i)$ for all
$i\in[n]$. A signed permutation $w$ of $[\pm n]$ acts on $\mb
R^n$ by permuting and negating coordinates:
\begin{align*}
w(v_1,v_2,\ldots,v_n) = \big(v_{w\inv(1)}, v_{w\inv(2)}, \ldots,
v_{w\inv(n)}\big), \text{ where } v_{-i} = -v_i \text{ for } i \in [n].
\end{align*}

The Coxeter group $D_n$ is the subgroup of the group of signed
permutations of $[\pm n]$ that negate an even number of elements of
$[n]$. The reflection arrangement $\Arrangement$ of $D_n$ consists of
the hyperplanes 
\begin{align*}
H_{ij} &= \{ \vec v \in \mb R^n : v_i = v_j \},
\text{ where } i\neq j,-i \in [\pm n].
\end{align*}

\subsection{Intersection lattice}
\label{ss:IntLatticeD}

A \defn{set partition} of $[\pm n]$ is a collection of nonempty
subsets $B = \{B_1, \ldots, B_r\}$ of $[\pm n]$ such that $\bigcup_i
B_i = [\pm n]$ and $B_i \cap B_j = \emptyset$ for $i \neq j$. The sets
$B_i$ in $B$ are called the \defn{blocks} of $B$. The collection of
set partitions of a finite set form a finite lattice, with partial
order given by: $B \leq C$ if and only if every block of $C$ is
contained in a block of $B$. If $A\subseteq [\pm n]$, then $\overline
A = \{ -a : a \in A \}$.

The intersection lattice $\IntersectionLattice$ of $\Arrangement$ is
isomorphic to the sublattice of set partitions of $[\pm n]$ of the
form $\{B_1, \ldots, B_r, C, \overline B_r, \ldots, \overline B_1\}$,
where $C$ satisfies: $C=\emptyset$ or $|C|\geq4$; and $C = \overline
C$ \cite[Theorem 4.1]{BarceloIhrig1999}. The isomorphism is given by
\begin{align*}
P = \{P_1, &\ldots, P_r\} \mapsto \\[-2ex]
&\Big\{\vec v \in V : v_i = v_j\text{ if } i,j \in P_h 
 \text{ for some }h\in[r] \Big\}
= \bigcap_{h=1}^{r} \left(\bigcap_{i,j \in P_h} H_{ij}\right),
\end{align*}
where $P$ is a set partition of $[\pm n]$ and $v_{-i} = -v_i$ for
$i\in\mb N$. 

To simplify notation, we let $\pi(X)$ denote the set partition of
$[\pm n]$ induced by $X\in\IntersectionLattice$, and we let
$\{B_1,\ldots,B_r;C\}$ denote the set partition $\{B_1, \ldots, B_r,
C, \overline B_r, \ldots, \overline B_1\}$. The block $C$ is called
the \defn{central block}. Under this isomorphism the action of $D_n$
on $X \in \IntersectionLattice$ is given by permuting the elements of
$\pi(X)$. That is, $\pi(w(X)) = w(\pi(X))$ for all $w\in D_n$ and
$X\in\IntersectionLattice$.

\subsection{Canonical basis}
\label{ss:CanonicalBasis}
The set partition $\pi(X) = \SetPartitionD$ describes a basis of
$X$. For each $i\in[r]$, let 
\begin{align*}
\bs\beta_i = \sum_{j\in B_i} \vec e_j,
\end{align*}
where $\vec e_1,\ldots,\vec e_n$ is the standard basis of $\mb R^n$
and $\vec e_{-j} = -\vec e_j$ for $j\in[n]$. The vectors
$\bs\beta_1,\ldots,\bs\beta_r$ form a basis of the subspace $X$
called the \defn{canonical basis} of $X$.

\subsection{The length of the longest path in $\InvariantQuiverD[2m]$}

This serves as a quick example to illustrate the approach we take in
the following section. 

The Coxeter group $D_{2m}$ contains an element $w_0$ that acts on $V$
by central reflection. That is, $w_0(\vec v) = -\vec v$ for all $\vec
v \in V$. Therefore, $\sigma_X(w_0) = (-1)^{\dim(X)}$. So if $A$ is an
arrow in $\Quiver$, then $w_0(A) = -A$. It follows from Lemma
\ref{l:ConditionForNoArrow} that there is no arrow $\Orbit \arrow
\Orbit'$ in $\InvariantQuiver[D_{2m}]$ if $\Orbit \lessdot \Orbit'$.
Combined with Proposition \ref{p:NoArrowsFromTop}, this implies that
the length of the longest path in $\InvariantQuiver[D_{2m}]$ is at
most $\frac{2m-1}2$, since $2m$ is the length of the longest path in
$\Quiver$. This establishes the following.

\begin{Proposition}
The length of the longest path in $\InvariantQuiverD[2m]$ is at most
$m-1$.
\end{Proposition}

This implies that the Loewy length of the descent algebra
$\DescentAlgebraD[2m]$, for $m\geq2$, is at most $m$ (see the proof
Theorem \ref{t:LoewyLengthTypeD}). Also, the same argument also gives
an upper bound of $\left\lceil \frac{n}2 \right\rceil$ for the Loewy
length of the descent algebra $\DescentAlgebra[B_n]$. These are both
equalities \cite[\S5E]{BonnafePfeiffer2006}. 

\subsection{The Loewy length of $\DescentAlgebraD$}

In this section we develop necessary conditions on $\Orbit_X$ and
$\Orbit_Y$ for there to be an arrow from $\Orbit_X$ to $\Orbit_Y$ in
$\InvariantQuiverD$.

For $X\in\IntersectionLattice$, let $\pi(X) = \SetPartitionD$ denote
the set partition induced by $X$ (see \S\ref{ss:IntLatticeD}). Let
$\Even X$ denote the number of $i\in[r]$ with $|B_i|$ even, and let
$\Odd X$ denote the number of $j\in[r]$ with $|B_j|$ odd.

\begin{Lemma} 
\label{l:InfoAboutEven}
If there is an arrow $\Orbit'\arrow\Orbit$ in $\InvariantQuiverD$,
then $\Even Y \leq \Even X$ for all $X\in\Orbit'$ and $Y\in\Orbit$.
\end{Lemma} 

\begin{proof}
We prove that if $\Even X < \Even Y$, then there is no arrow
$\Orbit_X\arrow\Orbit_Y$.

Suppose $P$ is a path in $\Quiver$ beginning at a vertex $X'$ in
$\Orbit_X$ and ending at a vertex $Y'$ in $\Orbit_Y$. Since $\pi(X')$
and $\pi(X)$ are in the same orbit, $\Even{X'} = \Even X$. Similarly,
$\Even{Y'} = \Even{Y}$. So $\Even{X'} < \Even{Y'}$. 

If every even-sized non-central block $B_i$ in $\pi(Y') =
\SetPartitionD$ contains an even-sized non-central block of $\pi(X')$,
then $\Even{X'} \geq \Even{Y'}$, contrary to our assumption.
Therefore, for some $i\in[r]$ the block $B_i\in\pi(Y')$ is even-sized
and is a union of an even number of odd-sized blocks of $\pi(X')$. 

Let $w$ be the signed permutation that negates all elements of $B_i$.
Then $w\in D_{2m+1}$ since $|B_i|$ is even and $w(Z) = Z$ for all
vertices of $P$. Since $w$ negates an even number of blocks of
$\pi(X')$ and fixes the others, $w$ negates an even number of vectors
(and fixes the others) in a canonical basis of $X'$. Hence
$\sigma_{X'}(w)=1$. Similarly, $w$ negates exactly 1 block of
$\pi(Y')$ while fixing the others, so $\sigma_{Y'}(w)=-1$. Hence,
$w(P)=-P$. By Lemma \ref{l:ConditionForNoArrow}
there is no arrow $\Orbit_X \arrow\Orbit_Y$ in $\InvariantQuiverD$.
\end{proof}

\begin{Lemma}
\label{l:ArrowImplies1OddBlock}
If $X'\arrow Y'$ is an arrow in $\Quiver$ and  if $\Odd{X'} \neq 1$, then
there is no arrow $\Orbit_{X'} \arrow \Orbit_{Y'}$ in $\InvariantQuiverD$. 
\end{Lemma}

\begin{proof}
Suppose $X\arrow Y$ is an arrow in $\Quiver$ with $X\in\Orbit_{X'}$ and
$Y\in\Orbit_{Y'}$. Since $\pi(X')$ is in the orbit of $\pi(X)$, $\Odd
X = \Odd{X'} \neq 1$. Since $Y\lessdot X$, $\pi(Y)$ is obtained from
$\pi(X)=\SetPartitionD$ by merging two blocks. Let
$\bs\beta_1,\ldots,\bs\beta_r$ be the canonical basis for $X$ (see
\S\ref{ss:CanonicalBasis}).

\medskip\noindent
\emph{Case 1.} $\pi(Y)$ is obtained from $\pi(X)$ by merging
$B_i$ and $B_j$, where $i\neq j$.
\smallskip

If $|B_i \cup B_j|$ is even, then let $w$ be the signed permutation
that negates the elements of $B_i$ and $B_j$. Then $w$ negates two
elements of the canonical basis of $X$ and fixes the other basis
elements, so $w(X)=X$ and $\sigma_{X}(w)=1$. Since
$\{\bs\beta_{i}+\bs\beta_{j}\} \cup \{\bs\beta_h : h \neq i,j \}$ is a
basis of $Y$, and since $w$ negates $\bs\beta_{i}+\bs\beta_{j}$ and
fixes the others, it follows that $w(Y)=Y$ and $\sigma_Y(w)=-1$.
Therefore, $w(X\arrow Y) = -(X\arrow Y)$.

If $|B_i \cup B_j|$ is odd, then let $w$ be the signed permutation
that negates the elements of $B_i,B_j$ and the elements of $B_h$,
where $h\neq i,j$ and $|B_h|$ is odd (such a block exists since $\Odd
X \neq 1$). Then $w$ negates three elements of the canonical basis of
$X$ and two elements of the basis
$\{\bs\beta_{i}+\bs\beta_{j}\}\cup\{\bs\beta_a:a\neq i,j\}$ of $Y$.
Therefore, $\sigma_X(w)=-1$ and $\sigma_Y(w)=1$. It follows that
$w(X\arrow Y)=-(X\arrow Y)$. 

\medskip\noindent
\emph{Case 2.} $\pi(Y)$ is obtained from $\pi(X)$ by merging
$B_i$ and $\overline B_j$, where $i\neq j$.
\smallskip

This is argued as is \emph{Case 1.}

\medskip\noindent
\emph{Case 3.} $\pi(Y)$ is obtained from $\pi(X)$ by merging
the blocks $B_i, \overline B_i$ and $C$. 
\smallskip

If $|B_i|$ is even, then let $w$ be the signed permutation that
negates the elements of $B_i$. Then $w$ negates one element of the
canonical basis of $X$ and no elements of the basis
$\{\bs\beta_a:a\neq i\}$ of $Y$. Thus, $w(X\arrow Y) = -(X\arrow Y)$.

If $|B_i|$ is odd, then let $w$ be the signed permutation that negates
the elements of $B_i$ and the elements of $B_h$, where $h\neq i$ and
$|B_h|$ is odd. Then $w$ negates two elements of the canonical basis
of $X$ and one element of the basis $\{\bs\beta_a:a\neq i\}$ of $Y$.
It follows that $w(X\arrow Y)=-(X\arrow Y)$. 

\medskip
It follows from Lemma \ref{l:ConditionForNoArrow} that there is no
arrow $\Orbit_{X'} \arrow \Orbit_{Y'}$ in $\InvariantQuiverD$.
\end{proof}



\begin{Proposition}
The length of the longest path in $\InvariantQuiverD$ is at
most $m+1$.
\end{Proposition}

\begin{proof}
Suppose $\Path[2] \Orbit l$ is a path in
$\InvariantQuiverD$. Then for $0\leq i\leq l$, there is $X_i \in
\Orbit_i$ such that $X_l \leq \cdots \leq X_1 \leq X_0$. 
Note that $X_0 \neq V$ by Proposition \ref{p:NoArrowsFromTop}. 

For each $j \in [l]$, let $d_j = \dimDiff[j]$. If $d_i \geq 2$ for all
$i\in[l]$, then
\begin{align*}
2l\leq\sum_{i=1}^l d_i=\dim(X_0)-\dim(X_l)\leq\dim(X_0)\leq2m,
\end{align*}
so $l\leq m$. Suppose instead that 
$d_j=1$ for some $j\in[l]$, and let $i$ be the smallest such
$j$. By the choice of $i$, $X_{i-1}\arrow X_i$ is an arrow in $\Quiver$
with $X_{i-1}\in\Orbit_{i-1}$ and $X_i\in\Orbit_i$. Then
$\Odd{X_{i-1}}=1$ by Lemma \ref{l:ArrowImplies1OddBlock} 
and $\Even{X_{i-1}} \leq \Even{X_0}$ by Lemma \ref{l:InfoAboutEven}. 

Recall that for each $X\in\IntersectionLattice$, if
$\pi(X)=\SetPartitionD$, then $\dim(X)=r$ (\S\ref{ss:CanonicalBasis}). In
particular, $\dim(X) = \Even{X}+\Odd{X}$ and $\dim(X)\leq
(2m+1)-\Even{X}$. Therefore, since 
$\Even{X_{i-1}} \leq \Even{X_0}$, 
\begin{align*}
\dim(X_0) 
\leq (2m+1)-\Even{X_0} 
\leq (2m+1)-\Even{X_{i-1}}.
\end{align*}
By the choice of $i$, $d_j \geq 2$ for all $j \in [i-1]$, so 
\begin{align*}
2(i-1) & \leq \dim(X_0)-\dim(X_{i-1})\\
&\leq 
\Big((2m+1)-\Even{X_{i-1}}\Big)
-\Big(\!\Even{X_{i-1}}+\Odd{X_{i-1}}\Big)\\
&\leq
2\big(m - \Even{X_{i-1}}\big).
\end{align*}
Since the length of $(\Orbit_{i-1}\arrows\Orbit_l)$ is bounded by
$\dim(X_{i-1})$, 
\begin{align*}
l & = (l - (i-1)) + (i-1) \\
& \leq 
\dim(X_{i-1}) + 
\big(m-\Even{X_{i-1}}\big) \\
& \leq
\big(\!\Even{X_{i-1}} + \Odd{X_{i-1}}\big) +
\big(m-\Even{X_{i-1}}\big) \\
& 
\leq m+1. \qedhere
\end{align*}
\end{proof}

\begin{Theorem}
\label{t:LoewyLengthTypeD}
For all $m\geq2$, the Loewy length of $\DescentAlgebraD$ is $m+2$.
\end{Theorem}
\begin{proof}
By Theorem \ref{t:Bidigare}, the Loewy length of $\DescentAlgebraD$ is
the Loewy length of $\InvariantSubalgebraD$. 
By Lemma \ref{l:LoewyLength} and the previous Proposition, the Loewy
length of $\InvariantSubalgebraD$ is bounded by $m+2$. By Corollary
5.9(b) of \cite{BonnafePfeiffer2006}, this is also a lower bound.
\end{proof}

\bibliographystyle{../../bibstyles/hapalike.bst}
\bibliography{references}

\end{document}